\documentclass[12pt]{article}

\usepackage{amsmath}

\usepackage{graphicx}
\usepackage{endfloat}
\usepackage{amssymb}
\usepackage{multirow}

\usepackage[default]{jasa_harvard}    
\usepackage{JASA_manu}

\DeclareMathOperator*{\argmin}{arg\,min}
\DeclareMathOperator*{\argmax}{arg\,max}
\def\x{{\bf x}}

\def\X{{\bf X}}
\def\btheta{{\boldsymbol \theta}}

\def\bgamma{{\boldsymbol \gamma}}
\def\bxi{{\boldsymbol \xi}}
\def\bupsilon{{\boldsymbol \upsilon}}
\def\bvartheta{{\boldsymbol \vartheta}}
\def\a{{\bf a}}
\def\1{{\bf 1}}
\def\g{{\bf g}}
\def\h{{\bf h}}
\def\A{{\bf A}}
\def\d{{\bf d}}
\def\e{{\bf e}}
\def\r{{\bf r}}

\newtheorem{thm}{Theorem}

\newtheorem{alg}{Algorithm}
\newtheorem{prop}{Proposition}

\newcommand{\tr}{^{\tt T}}

\begin{document}


\title{Bayes Multiple Decision Functions{\normalsize {\Large }}}
\author{Wensong Wu
 \\Department of Mathematics and Statistics\\Florida International University, Miami, FL, 33199\\\texttt{wensong.wu@fiu.edu}
\and
Edsel A. Pe\~na\\Department of Statistics\\University of South Carolina, Columbia, SC 29208\\\texttt{pena@stat.sc.edu}}
\date{\today}

\maketitle



\newpage
\begin{center}
\textbf{Abstract}
\end{center}
This paper deals with the problem of simultaneously making many ($M$) binary decisions
based on one realization of a random data matrix $\mathbf{X}$. $M$ is typically large and
$\mathbf{X}$ will usually have $M$ rows associated with each of the $M$ decisions to make, but
for each row the data may be low dimensional. Such problems arise in many practical arenas
such as in the biological and medical sciences, where the available dataset is from microarrays or
other high-throughput technology and with the goal being to decide which among of many genes are relevant with respect to some phenotype of interest; in the engineering and reliability
sciences; in astronomy; in education; and in business. A Bayesian decision-theoretic
approach for this problem is implemented with the overall loss function being a cost-weighted linear combination
of Type I and Type II loss functions. The class of loss functions considered
allows for the use of the false discovery rate (FDR), false nondiscovery rate (FNR), and
missed discovery rate (MDR) in assessing the decision. Through
this Bayesian paradigm, the Bayes multiple decision function (BMDF) is derived and an efficient
algorithm to obtain the optimal Bayes action is described. In contrast to many works in the
literature where the rows of the matrix $\mathbf{X}$ are assumed to be stochastically independent,
we allow in this paper a dependent data structure with the associations obtained through a class of
frailty-induced Archimedean copulas. In particular, non-Gaussian dependent data structure,
which is the norm rather than the exception when dealing with failure-time data, can be
entertained. The numerical implementation of the determination of the Bayes optimal action
is facilitated through sequential Monte Carlo techniques. The main theory developed could also
be extended to the problem of multiple hypotheses testing, multiple classification and
prediction, and high-dimensional variable selection. The proposed procedure is illustrated
for the simple versus simple and for the composite hypotheses setting
via simulation studies. The procedure is also applied to a subset of a real microarray data
set from a colon cancer study.


\noindent\textsc{\bf Keywords}: {Archimedean Copula, Bayes framework, Decision-theoretic framework, False discovery proportion, False Discovery Rate, Frailty, Multiple testing,  Sequential Monte Carlo.}

\newpage

\section{Introduction}\label{sec: intro}

The advent of computer-automated high-throughput data-gathering technology, epitomized by the microarray,
has led to the generation of so-called ``large $M$, small $n$'' data sets, which are those characterized by a
large number, $M$, of variables (hereon called {\em genes} for historical reasons), 
which are observed or measured on a relatively small
number, $n$, of subjects or units. Examples of such data sets in different scientific fields could, for instance, be
found in  \citeasnoun{HasTibFri09} and \citeasnoun{Efr08}. 

For such data sets a typical goal is to choose an action, $a_m$, associated with gene $m$,
from a set of possible actions, $\mathcal{A}_m$, for each of the $M$ genes.
For example, in a two-group microarray data set, one may want to decide, for each gene, whether
it is differentially expressed between the two groups (action is $a = 1$), or whether it is not differentially
expressed between the two groups (action is $a = 0$). This situation corresponds to the problem of simultaneously testing
multiple pairs of null and alternative hypotheses. 

In this paper we shall focus on these two-point action spaces for each of the genes,
that is, those with $\mathcal{A}_m = \{0,1\}$.
Of interest therefore is to choose a vector of actions
\begin{displaymath}
\mathbf{a} = (a_1,a_2, \ldots,a_M)\tr \in \mathcal{A} = \{0,1\}^M
\end{displaymath}
based on the observed ``large $M$, small $n$'' data set. For the $m$th gene there will be associated a $\theta_m \in \{0,1\}$, which is unknown, 
representing the {\em correct} action to take, which unfortunately is unknown. Thus, for the $M$ genes there will be an 
unknown vector
\begin{displaymath}
\btheta = (\theta_1, \theta_2, \ldots, \theta_M)\tr \in \Theta = \{0,1\}^M
\end{displaymath}
representing the vector of actions that {\em ought} to be taken. This ${\btheta}$ will be referred to as the 
{\em state of reality}. In light of this state of reality vector $\btheta$, a chosen action vector $\mathbf{a}$ 
will have consequences quantified through a loss. That is, there will be a mapping
\begin{displaymath}
(\mathbf{a},{\btheta}) \mapsto L(\mathbf{a},{\btheta})
\end{displaymath}
where $L(\mathbf{a},\mathbf{\theta})$ is the loss that is incurred with the action $\mathbf{a}$ when reality is
${\btheta}$. Such a loss must take into account the loss incurred when the action is $a = 1$ when
reality is $\theta = 0$, called a Type I error, as well as the loss incurred when the action is $a = 0$ when
reality is $\theta = 1$, called a Type II error. There will be a variety of ways of measuring 
the overall Type I and Type II errors in such multiple decision problem, which will be formally described in Section 2.

These multiple decision problems appertaining to such ``large $M$, small $n$'' data sets therefore lend naturally to a
decision-theoretic framework which will be discussed in more detail in Section 2. 
In addition to this decision-theoretic framework, we implement a Bayesian approach to decision-making by putting a
prior probability on the unknown state of reality ${\btheta}$. Coupled with the appropriate loss function, we
obtain the Bayes multiple decision action. To achieve this, we derive the mathematical form of the Bayes multiple
decision function, abbreviated BMDF, and describe an efficient computational implementation of this BMDF under
varied scenarios in terms of loss functions and data structures.

A decision-theoretic and a Bayesian approach to these multiple decision problems with high-dimensional data is certainly not
new as can be seen from the papers by \citeasnoun{Mul04}, \citeasnoun{ScoBer06} and \citeasnoun{Sar08}. Other approaches include, but not limited to, \citeasnoun{Sto03}, \citeasnoun{SunCai07} and \citeasnoun{PenHabWu10}. See also the monograph \citeasnoun{Efr10mono}. An innovative major
contribution of this paper is the use of a general class of loss functions that encompasses many of the loss functions
that have been used in earlier works. For instance, the general class of loss functions introduced in Section 2 includes
as special cases those that involve false positives and false negatives as well as the commonly-used false discovery rates 
and false nondiscovery rates. Another major contribution is an efficient algorithm for computationally finding the
Bayes multiple decision action, an algorithm that has computational order of at most $O(M^2\log M)$. Many papers have 
dealt with the situation where the observables from each of the genes are stochastically independent. We go beyond this
usual assumption by incorporating dependencies among these observables, with the dependence structure induced by 
frailty-type models, which also takes the form of Archimedean
copulas. This dependent modeling approach utilizes ideas from the survival analysis area.
The statistical models
governing these multiple decision problems possess more complications than the description above. This is so since, even though
the parameter of main interest is the state of reality vector $\btheta$, there will be other unknown model parameters that are
present which may be viewed as nuisance, and we need to be able to deal with these nuisance parameters as well, both in
constructing the Bayes multiple decision functions and in the prior probability specification. In our development of the BMDF
we therefore first consider the situation of a simple null hypothesis versus simple alternative hypothesis setting wherein the distributional model for the random
observable for the $m$th gene is completely known under either $\theta_m = 0$ or $\theta_m = 1$; and then utilize the
results for this setting to solve the problem for the composite versus composite setting where unknown nuisance
parameters are present. An interesting development is the use of Sequential Monte Carlo (SMC) techniques to numerically approximate the Bayes multiple decision action especially in the presence
of dependencies among the observables and also in the prior probability specification.

We now outline the contents of this paper. In Section \ref{sec: setting} we will introduce the mathematical setting and elements 
of the multiple decision problem, including a general class of loss functions. Section \ref{sec: BMDF} will demonstrate the general 
form of the BMDF along with a computationally efficient algorithm of finding this BMDF in both simple and composite hypotheses 
testing settings. Section \ref{sec: concrete} will give the expressions of the BMDF under three concrete loss functions. We will 
introduce the frailty-based dependent data and models in Section \ref{sec: data dep}, and in Section \ref{sec: computation} 
we will discuss the computational aspects of the posterior calculations under the dependent models, and give some algorithms 
using Sequential Importance Sampling (SIS). In Section \ref{sec: illustration} we will illustrate the BMDF in 
some concrete multiple decision problems, and compare the performance with currently used procedures via simulation studies. We will also apply the BMDF to a subset of a microarray data set.
We will conclude the paper in Section \ref{sec: conclude} with some remarks.

\section{Decision-Theoretic Elements}\label{sec: setting}
\subsection{Multiple Decision Problem}\label{general problem}
Let $(\Omega,\mathcal{F},\mathcal{P})$ be a basic statistical model with $\mathcal{P}$ being a collection of probability measures. For $m=1,2,\dots, M$, where $M\ge 1$ is a known integer, let $X_m:(\Omega,\mathcal{F})\to(\mathcal{X}_m,\mathcal{B}_m)$, where $\mathcal{X}_m$ is some space and $\mathcal{B}_m$ is an associated $\sigma$-field of subset, of $\mathcal{X}_m$. Let $\X=(X_1,X_2,\dots,X_M):(\Omega,\mathcal{F})\to(\mathcal{X},\mathcal{B})$, where $\mathcal{X}=\bigotimes_{m=1}^M \mathcal{X}_m$ is the product sample space and $\mathcal{B}$ is the associated product $\sigma$-field. A realization $\X=\x$ will be called a (sample) data. For any $P\in\mathcal{P}$, the induced joint probability measure of $\X$ is $Q=P\X^{-1}$, whereas the marginal probability measure of $X_m$ is $Q_m=PX_m^{-1}$. Let $\mathcal{Q}=\{P\X^{-1}:P\in\mathcal{P}\}$ denote the collection of all probability measures of $\X$. Consider a mapping $\bvartheta=(\vartheta_1,\vartheta_2,\dots,\vartheta_M)\tr :\mathcal{Q}\to\{0,1\}^M$, where $\vartheta_m:\mathcal{Q}\to\{0,1\}$ only depends on the $m$th marginal probability measure. In essence, the parameter of main interest is $\bvartheta(Q)=\btheta=(\theta_1,\theta_2,\dots,\theta_M)\tr$, which takes values in the parameter space $\Theta=\{0,1\}^M$. Observe that $\mathcal{Q}$ can be decomposed via $\mathcal{Q}=\biguplus_{\btheta\in\Theta} \mathcal{Q}_{\btheta}$, where $\mathcal{Q}_\btheta=\{Q\in\mathcal{Q}:\bvartheta(Q)=\btheta\}, \forall\btheta\in\Theta$. Let $\a=(a_1,a_2,\dots,a_M)\tr$ be an action in the action space $\mathcal{A}=\{0,1\}^M$. Let $L^\circ:\mathcal{A}\times\mathcal{Q}\to\mathbb{R}$ be a loss function such that 
$L^\circ(\a,Q)=L(\a,\bvartheta(Q)), \forall\a\in\mathcal{A}, \forall Q\in\mathcal{Q},$
where $L:\mathcal{A}\times\Theta\to\mathbb{R}$ is a loss function belonging to a class $\mathcal{L}$ which will be discussed in Section \ref{sec:loss function}.

Let $\mathcal{P}(\mathcal{A})$ be the collection of all probability measures over $\mathcal{A}$.
A randomized decision function is a measurable function 
$\delta:\mathcal{X}\to\mathcal{P}(\mathcal{A})$ 
such that $\delta(\x)(\a)\equiv \delta(\x)(\{\a\})$ denotes the probability of choosing action $\a\in\mathcal{A}$ given data $\x$. Therefore, to make a specific decision via $\delta$ given data $\x$, an element $\a\in\mathcal{A}$ will be chosen according to the probability mass function $\delta(\x)(\cdot)$. Let $\mathcal{D}$ denote the space of all such randomized decision functions. Note that a nonrandomized decision function $\delta ^*:\mathcal{X}\to\mathcal{A}$, which maps $\x\in\mathcal{X}$ directly to an $\a\in\mathcal{A}$ belongs to $\mathcal{D}$, since it is just a degenerate probability measure over $\mathcal{A}$, that is, $\delta(\x)(A)=I(\delta(\x)\in A)$. 

The risk of a randomized decision function $\delta\in\mathcal{D}$ for any $Q\in\mathcal{Q}$ associated with a loss $L^\circ$ is defined by
\begin{equation*}
R^\circ(\delta,Q)=E_{\X\sim Q}\, E_{\A\sim \delta(\X)(\cdot)}\,L^\circ(\A,Q)= \sum_{\btheta\in\Theta}\,\left[\int_{\mathcal{X}}\sum_{\a\in\mathcal{A}}\,L(\a,\btheta)[\delta(\x)(\a)]Q({\rm{d}}\x)\right]\,I(\bvartheta(Q)=\btheta).
\end{equation*}

Now, assume that for any $\btheta\in\Theta$, $\mathcal{Q}_\btheta$ is an identifiable parametric class given by $\mathcal{Q}_\btheta=\{Q_\btheta(\cdot \, ;\,\gamma_\btheta):\gamma_\btheta\in\Gamma_\btheta\}$, where $\gamma_\btheta$ is a nuisance parameter. This implies that $\mathcal{Q}$ is an identifiable model with respect to the parameter $(\btheta,\gamma_\btheta)$ which belongs to the enlarged parameter space $\Theta^\circ=\cup_{\btheta\in\Theta}[\{\btheta\}\times\Gamma_\btheta]$. Then, for any $\delta\in\mathcal{D}$ and $Q\in\mathcal{Q}$, the risk function is given by
$R^\circ(\delta,Q)=\sum_{\btheta\in\Theta}R(\delta,(\btheta,\gamma_\btheta))I(\bvartheta(Q)=\btheta, Q=Q_\btheta(\cdot,\,;\,\gamma_\btheta)),$
where $R:\mathcal{D}\times\Theta^\circ\to\mathbb{R}$ is given by 
$$R(\delta,(\btheta,\gamma_\btheta))=E_{\X\sim Q_\btheta(\cdot;\gamma_\btheta)}\, E_{\A\sim \delta(\X)(\cdot)}\,L(\A,\btheta)=\int_{\mathcal{X}}\sum_{\a\in\mathcal{A}}\,L(\a,\btheta)[\delta(\x)(\a)]Q_\btheta({\rm{d}}\x;\gamma_\btheta).$$

Furthermore, the decomposition of $\mathcal{Q}$ becomes $\mathcal{Q}=\biguplus_{\btheta\in\Theta}\Big\{Q_\btheta(\cdot;\gamma_\btheta):\bvartheta(Q_\btheta(\cdot;\gamma_\btheta))=\btheta,\gamma_\btheta\in\Gamma_\btheta\Big\}.$ A prior probability measure on $\mathcal{Q}$ can be constructed by specifying a prior probability measure on $(\Theta^\circ,\sigma(\Theta^\circ))$, where $\sigma(\Theta^\circ)$ is the $\sigma$-field generated by a semi-ring $\mathcal{C}=\{\{\btheta\} \times C_{\btheta}: C_{\btheta}\in\sigma(\Gamma_\btheta), \btheta\in\Theta\}$ with $\sigma(\Gamma_\btheta)$ being the $\sigma$-field on $\Gamma_\btheta$. Define a probability measure $\Pi^*$ on $\mathcal{C}$ such that for any $\btheta_0\in\Theta$ and any $C_{\btheta_0}\in\sigma(\Gamma_{\btheta_0})$,
$$\Pi^*(\btheta=\btheta_0,\gamma_\btheta\in C_{\btheta_0})=\Pi(\btheta_0)\int_{C_{\btheta_0}} P_{\btheta_0}({\rm{d}}\gamma_{\btheta_0})=\Pi(\btheta_0)P_{\btheta_0}(C_{\btheta_0}),$$ where  $\Pi(\cdot)$ is a probability measure on $\Theta$ and $P_\btheta(\cdot)$ is a probability measure on $\Gamma_\btheta$. This induces, by Caratheodory's extension theorem, a prior probability measure $\Pi^\circ$ on $(\Theta^\circ, \sigma(\mathcal{C})=\sigma(\Theta^\circ))$. Since the elements of $\mathcal{Q}$ are identified by $(\btheta,\gamma_\btheta)$, there is a one-to-one mapping $h: \Theta^\circ\to\mathcal{Q}$ with $(\btheta,\gamma_\btheta) \stackrel{h}{\leftrightarrow}
 Q_\btheta(\cdot;\gamma_\btheta)$. Therefore, $\Pi^\circ$ determines a prior measure on $(\mathcal{Q}, \sigma(\mathcal{Q}))$, where $\sigma(\mathcal{Q})=h\,\sigma(\Theta^\circ).$ The Bayes risk function of a randomized decision function $\delta\in\mathcal{D}$ for a prior $\Pi^\circ$ is defined via
\begin{equation*}
r_{\Pi^\circ}(\delta)=E_{Q\sim \Pi^\circ}\,R^\circ(\delta,Q)=\int_\Theta\int_{\Gamma_{\btheta} }R(\delta,(\btheta,\gamma_\btheta)) P_\btheta ({\rm{d}}\gamma_\btheta)\Pi({\rm{d}}\btheta).
\end{equation*}
A randomized decision function $\delta^*$ is called a Bayes multiple decision function (BMDF) if
\begin{equation*}
\delta^*=\argmin_{\delta\in\mathcal{D}}\, r_{\Pi^\circ}(\delta).
\end{equation*}
The multiple decision problem is to find the BMDF, which is the optimal procedure for choosing the $M$-dimensional action vector from the Bayesian decision theory perspective. More practically, there is the issue of finding the optimal action in a computationally efficient manner.

\subsection{Loss Functions}\label{sec:loss function} 
The loss function $L: \mathcal{A}\times \Theta\to\mathbb{R}$
quantifies the error committed for a pair $(\a,\btheta)$ of action and state of reality.  We shall consider a class of cost-weighted loss functions with elements
\begin{equation}\label{general loss}
L(\a,\btheta)=C_0 L_0(\a,\btheta)+C_1 L_1(\a,\btheta),
\end{equation}
where $C_0\ge 0$ and $C_1\ge 0$ are pre-determined costs for loss functions $L_0$ and $L_1$, respectively. 
The generic forms of $L_0$ and $L_1$ are
$$L_0(\a,\btheta)=[\alpha_0(\a\tr\1)\g_0(\a)]\tr\,[\beta_0(\btheta\tr\1)\h_0(\btheta)],$$
$$L_1(\a,\btheta)=[\alpha_1(\a\tr\1)\g_1(\a)]\tr\,[\beta_1(\btheta\tr\1)\h_1(\btheta)],$$
where $\alpha_j:\mathbb{R}\to\mathbb{R}$, $\beta_j:\mathbb{R}\to\mathbb{R}$, $\g_j:\mathcal{A}\to\mathcal{A}$, and $\h_j:\Theta\to\Theta$ for $j=0,1$. We assume further that, for $j=0,1$, $\g_j$ is $\tau_j$-invariant with respect to the sub-action space $\mathcal{A}_k\equiv \{\a\in\mathcal{A}:\a\tr\1=k\}$ for $k\in\mathcal{M}\equiv \{0,1,\dots,M\}$, in the sense that there exists a mapping $\tau_j:\mathcal{M}\to\mathcal{M}$ associated with $\g_j$ such that $\a\in\mathcal{A}_k$ implies $\g_j(\a)\in\mathcal{A}_{\tau_j(k)}$. Examples of $\tau$ include the identity mapping with $\tau_0(k)=k$ and also $\tau_1(k)=M-k$. Then $\g_0$ with $\g_0(\a)=\a$ is $\tau_0$-invariant, while $\g_1$ with $\g_1(\a)=\1-\a$ is $\tau_1$-invariant. With $a\vee b=\max(a,b)$, some examples of loss functions $L_0$ and $L_1$ on $\mathcal{A}\times\Theta$ are
the False Positive Proportion (FP): $(\a,\btheta)\mapsto\dfrac{\a\tr(\1-\btheta)}{M};$
the False Negative Proportion (FN) $(\a,\btheta)\mapsto\dfrac{(\1-\a)\tr\btheta}{M};$
the False Discovery Proportion (FDP) $(\a,\btheta)\mapsto\dfrac{\a\tr(\1-\btheta)}{(\a\tr\1)\vee 1};$
the Missed Discovery Proportion (MDP) $(\a,\btheta)\mapsto\dfrac{(\1-\a)\tr\btheta}{(\btheta\tr\1)\vee 1};$
and the False Nondiscovery Proportion (FNP) $(\a,\btheta)\mapsto\dfrac{(\1-\a)\tr\btheta}{((\1-\a)\tr\1)\vee 1}.$

Besides the aforementioned properties of $L_0$ and $L_1$, we also assume that $L_0$ and $L_1$ possess a complementarity property given by
$\g_1(\a)=\a_0-A_1\g_0(\a),\textrm{ for some }\a_0\in\mathcal{A} \textrm{ with } A_1>0.$
This property describes a relation between $L_0$ and $L_1$  which indicates that they are loss functions having complementary behaviors. For example, FP and FDP are proportions of false discoveries, where a discovery at the $m$th coordinate is having $a_m=1$, whereas FN, MDP, and FNP are proportions of false nondiscoveries. In the sequel, we will consider the pairs (FP, FN), (FDP, MDP), and (FDP, FNP) for ($L_0,L_1$) in the multiple decision problem. The cost constants $C_0$ and $C_1$ will generally be determined by the decision maker or subject matter specialist, and they reflect the consequences of false discoveries and false nondiscoveries.

\subsection{Multiple Testing Problems}\label{subsec: testing}
The simple-versus-simple multiple hypotheses testing problem is a particular case of this multiple decision problem. Suppose that the marginal probability distribution of $X_m$ satistisfies $Q_m\in\{Q_{m0},Q_{m1}\}$ with $Q_{m0}\neq Q_{m1}$, and the parameter vector is $\btheta=(I(Q_m=Q_{m1}),m=1,2,\dots,M)$. We may consider simultaneously the $M$ pairs of simple versus simple hypotheses $H_{m0}: Q_m=Q_{m0}$ versus $H_{m1}: Q_m=Q_{m1}$ for $m=1,2,\dots,M$. In this case, for $m=1,2,\dots,M$, $\theta_m=1(0)$ indicates whether $H_{m1}$ is (not) true, and the action $a_m=1(0)$ means rejecting (not rejecting) $H_{m0}$. 

Usually, independent Bernoulli priors are assigned to $\btheta$. Let $\pi_{m0},\pi_{m1}\in(0,1)$ be such that $\pi_{m0}+\pi_{m1}=1$ for $m=1,2,\dots,M$.  The prior probability mass function $\pi$ on $\btheta$ is specified by 
\begin{equation}\label{indep prior}
\pi(\btheta)=\prod_{m=1}^M \pi_{m0}^{1-\theta_m}\pi_{m1}^{\theta_m}I(\theta_m\in\{0,1\}).
\end{equation}
In this situation, the Bayes risk function of a randomized decision function $\delta$ for a prior mass function $\pi$ of $\btheta$ is 
$r_{\pi}(\delta)=E_{\btheta\sim \pi}\,E_{X\sim Q_\btheta}\, E_{\A\sim \delta(\X)(\cdot)}\, L(\A,\btheta),$
associated with a loss function $L\in\mathcal{L}$, where $Q_\btheta$ is the joint probability function of $\X$, given $\btheta$, and whose $m$th marginal distribution function is $Q_m=Q_{m\theta_m}$.

Suppose that $Q_m$, the marginal distribution of $X_m$, is in a class of distributions $\mathcal{Q}_m$ given by
$\mathcal{Q}_m=\{Q_m(\cdot\,;\gamma_m,\xi_m):\gamma_m\in\Gamma_m,\xi\in\Xi_m\}.$
Assume, for $m=1,2,\dots,M$, $\Gamma_m=\Gamma_{m0}\cup\Gamma_{m1}$ and $\Gamma_{m0}\cap\Gamma_{m1}=\varnothing$. Then $\mathcal{Q}_m$ has two subclasses denoted by
$$\mathcal{Q}_{m0}=\{Q_m(\cdot\,;\gamma_m,\xi_m):\gamma_{m}\in\Gamma_{m0},\xi\in\Xi_m\} \textrm{ and }
\mathcal{Q}_{m1}=\{Q_m(\cdot\,;\gamma_m,\xi_m):\gamma_{m}\in\Gamma_{m1},\xi\in\Xi_m\}.$$
Consider the $M$ pairs of composite hypotheses $H_{m0}:\,Q_m\in\mathcal{Q}_{m0}$ versus $H_{m1}:\,Q_m\in\mathcal{Q}_{m1}$, for $m=1,2,\dots,M$. Note that $\Gamma_{m0}$, $\Gamma_{m1}$, or $\Xi_m$ could be the same for all $m$, though in general they may be different. 
Let $\btheta=(I(Q_m\in\mathcal{Q}_{m1}),m=1,2,\dots,M)=(I(\gamma_m\in\Gamma_{m1}),m=1,2,\dots,M)\in\Theta=\{0,1\}^M$, $\bgamma=(\gamma_1,\gamma_2,\dots,\gamma_M)\in\Gamma\equiv\bigotimes_{m=1}^M \Gamma_m$, and $\bxi=(\xi_1,\xi_2,\dots,\xi_M)\in\Xi\equiv\bigotimes_{m=1}^M \Xi_m$. Also, for any $\btheta\in\Theta$, let $\Gamma_\btheta\equiv \bigotimes_{m=1}^M \Gamma_{m\theta_m}$. Then the extended parameter vector is
\begin{equation}\label{extended par vec}
(\btheta,\bgamma,\bxi)\in\Theta^\circ\equiv \biguplus_{\btheta\in\Theta}\{\{\btheta\}\times \Gamma_\btheta\times\Xi\},
\end{equation}
where $\btheta\in\Theta$ is the parameter of main interest in the multiple decision problem, while $(\bgamma,\bxi)\in\Gamma\times\Xi$ are nuisance parameters. Note that the value of $\bgamma$ determines the value of $\btheta$, but we only want to determine whether $\gamma_m\in\Gamma_{m0}$ or $\gamma_m\in\Gamma_{m1}$ rather than estimating the exact values of $\gamma_m$s. The parameter $\btheta\in\Theta$ and the action $\a\in\mathcal{A}$ have the same interpretation in this composite hypotheses testing problem as in the simple-vs-simple setting. 


Assume that the prior distribution on the enlarged parameter space $\Theta^\circ$ is 
\begin{equation}\label{composite prior}
\pi(\btheta,\bgamma,\bxi)=\prod_{m=1}^M(\pi_{m0}p_{m0}(\gamma_m,\xi_m))^{1-\theta_m}(\pi_{m1}p_{m1}(\gamma_m,\xi_m))^{\theta_m}.
\end{equation}  
where $p_{m0}$ and $p_{m1}$ are prior densities on $\Gamma_{m0}\times\Xi$ and $\Gamma_{m1}\times\Xi$, respectively, and with $\pi_{m0},\pi_{m1}\in\{0,1\}$ and $\pi_{m0}+\pi_{m1}=1$. 
The Bayes risk function of a randomized decision function $\delta$ for a prior density $\pi$ of $(\btheta,\bgamma,\bxi)$ associated with a loss function $L\in\mathcal{L}$ is 
$r_{\pi}(\delta)=E_{(\btheta,\bgamma,\bxi)\sim \pi}\,E_{\X\sim Q_\btheta(\cdot;\gamma,\bxi)}\, E_{\A\sim \delta(\X)(\cdot)}\, L(\A,\btheta),$
where $Q_\btheta(\cdot;\bgamma,\bxi)$ is the joint probability function of $\X$, given $(\btheta,\bgamma,\bxi)$, and the marginal probability measures are $Q_m=Q_m(\cdot;\gamma_{m},\xi_m)$, $m=1, 2, \dots, M$. To indicate that the marginal $Q_m\in\mathcal{Q}_{m\theta_m}$, we shall denote it by $Q_{m\theta_m}(\cdot;\gamma_m,\xi_m)$.

\section{Bayes Multiple Decision Functions}\label{sec: BMDF}
\subsection{BMDF in Simple Hypotheses}\label{subsec: BMDF simple}
Let $\pi(\cdot)$ be a prior probability mass function of $\btheta\in\Theta$. Then 
The Bayes risk function of $\delta$ for the prior $\pi$ is given by $r_\pi(\delta)=E_{\btheta\sim\pi}\,R(\delta,\btheta)=E_\btheta \,E_{\X|\btheta} \,E_{\A\sim \delta(\X)(\cdot)}\,L(\A,\btheta)=E_\X \,E_{\btheta|\X}\, E_{\A\sim \delta(\X)(\cdot)}\,L(\A,\btheta)=E_\X\, E_{\A\sim \delta(\X)(\cdot)}\, E_{\btheta|\X}\,L(\A,\btheta)$,
where $E_{\btheta|\X}$ is the expectation with respect to the posterior distribution of $\btheta$ given $\X$. 
For $\a\in\mathcal{A}$ and $\X=\x\in\mathcal{X}$, define the posterior expected loss by,
\begin{equation}\label{tilde L}
\tilde{L}(\a,\x)=E_{\btheta|\X=\x}\,L(\a,\btheta),
\end{equation}
and denote the optimal action when $\X=\x$ by
\begin{equation}\label{original a*}
\a^*(\x)=\argmin_{\a\in\mathcal{A}}\,\tilde{L}(\a,\x).
\end{equation}
Then the BMDF is
\begin{equation}\label{decision form}
\delta^*(\X)(B)=I(\a^*(\X)\in B),\quad B\in \sigma(\mathcal{A}).
\end{equation}

Notice that $\delta^*$ is degenerate at the nonrandomized decision function $\delta:\mathcal{X}\to \mathcal{A}$ with $\delta(\x)=\a^*(\x)$, which implies that we may always choose the BMDF to be nonrandomized. It is also worth noting that finding the optimal action $a^*$, and thus the BMDF $\delta^*$, via equation (\ref{original a*}) involves searching for the minimizer of the function $\tilde{L}$ among all $2^M$ elements of $\mathcal{A}$. When $M$ is relatively large, the searching order of $O(2^M)$ would become practically infeasible. Furthermore, note that this computational problem does not yet include the problem of computing the posterior distribution of $\btheta$ given $\X=\x$.

The idea for obtaining a computationally efficient algorithm to find the optimal action is to first find the restricted optimal action over the sub-action space $\mathcal{A}_k=\{\a\in\mathcal{A}:\a\tr\1=k\}$ for each $k\in\mathcal{M}$, and then to find the optimal action among these restricted optimal actions. Before presenting the results, we first define some relevant quantities that will be used. With the notation of the loss $L\in\mathcal{L}$ described in Section \ref{sec:loss function}, for $k\in\mathcal{M}$, let
\begin{align}
\d_0(k,\x)&=C_0\alpha_0(k)E_{\btheta|\X=\x}\,[\beta_0(\btheta\tr\1)\h_0(\btheta)]\label{func d0};\\
\d_1(k,\x)&=C_1\alpha_1(k)E_{\btheta|\X=\x}\,[\beta_1(\btheta\tr\1)\h_1(\btheta)]\label{func d1};\\
\e(k,\x)&=\d_0(k,\x)-A_1\d_1(k,\x)\label{func e};
\end{align}
and let $\r(k,\x)=(r_m(k,\x),\,m=1,2,\dots,M)$ be the rank vector of $\e(k,\x)$.
Also, define
\begin{equation}\label{func H}
H(k,\x)=\a_0\tr \d_1(k,\x)+\sum_{m=1}^M I(r_m(k,\x)\leq \tau_0(k))e_m(k,\x),
\end{equation}
where we recall that $\tau_0:\mathcal{M}\to\mathcal{M}$ is a mapping such that $\a\in\mathcal{A}_k$ implies $\g_0(\a)\in\mathcal{A}_{\tau_0(k)}$.

The following theorem describes a computationally efficient algorithm for finding the optimal action vector.

\begin{thm}\label{general solution}
For a multiple decision problem with loss function $L\in\mathcal{L}$ and prior probability mass function $\pi$ on $\btheta$, let $k^*:\mathcal{X}\to\mathcal{M}$ be defined via
\begin{equation*}
k^*(\x)=\argmin_{k\in\mathcal{M}}\,H(k,\x), \quad \x\in\mathcal{X},
\end{equation*}
where $H:\mathcal{M}\times\mathcal{X}\to\mathbb{R}$ is defined in (\ref{func H}). Then the BMDF is of the form (\ref{decision form}) with $a^*(\x)$ satisfying
\begin{equation}\label{Eqn Thm1}
\g_0(\a^*(\x))=\bigg(I\{r_m(k^*(\x),\x)\leq \tau_0(k^*(\x)\},\quad m=1,2,\dots,M\bigg).
\end{equation}
The searching order for obtaining this BMDF is no more than $O(M^2\log M)$.
\end{thm}

{\noindent\bf Proof: } Associated with the loss function $L$, $\tilde{L}$ defined in equation (\ref{tilde L}) has a specific form given by
$\tilde{L}(\a,\x)=\sum_{j=0}^1 C_j[\alpha_j(\a\tr\1)\g_j(\a)]\tr\,E_{\btheta|\X=\x}\,[\beta_j(\btheta\tr\1)\h_j(\btheta)].$
Restricting $\a$ on $\mathcal{A}_k$, 
\begin{align*}
\tilde{L}(\a,\x)&= \g_0(\a)\tr \d_0(k,\x)+\g_1(\a)\tr \d_1(k,\x)\\
                &=\a_0\tr \d_1(k,\x)+\g_0(\a)\tr[\d_0(k,\x)-A_1\d_1(k,\x)]\\
                &=\a_0\tr \d_1(k,\x)+\g_0(\a)\tr \e(k,\x),              
\end{align*}  
where $d_0(k,\x)$, $d_1(k,\x)$, and $\e(k,\x)$ are as defined in (\ref{func d0})-(\ref{func e}). 

Since for $\a\in\mathcal{A}_k$, $\g_0(\a)\tr\1=\tau_0(k)$, the optimal action on $\mathcal{A}_k$, denoted by $\a_k^*(\x)$, which minimizes $\tilde{L}(\a,\x)$ for $\a\in\mathcal{A}_k$, therefore satisfies
\begin{equation*}
\g_0(\a_k^*(\x))=\bigg(I\{r_1(k,\x)\leq \tau_0(k)\},I\{r_2(k,\x)\leq \tau_0(k)\},\dots,I\{r_M(k,\x)\leq \tau_0(k)\}\bigg)\tr
\end{equation*}
where we recall that, for $m=1,2,\dots,M$, $r_m(k,\x)$ is the rank of $e_m(k,\x)$ among the elements of $\e(k,\x)$, $m=1,2,\dots,M$. 
Thus 
$\tilde{L}(\a_k^*(\x),\x)=\a_0\tr \d_1(k,\x)+\sum_{m=1}^M I(r_m(k,\x)\leq \tau_0(k))e_m(k,\x),$
which equals  the function $H(k,\x)$. Therefore, for the $k^*(\x)$ in the statement of Theorem \ref{general solution}, $\a^*_{k^*(\x)} (\x)$ minimizes $\tilde{L}(\a,\x)$ over all actions $\a\in\mathcal{A}$. The optimal action, given $\X=\x$, is therefore
$\a^*(\x)=\a^*_{k^*(\x)} (\x)$,
which satisfies (\ref{Eqn Thm1}).

To see the computational order of the algorithm, for $k\in\mathcal{M}$, observe that to find $\a^*_k(\x)$, it is only necessary to know which are the $\tau_0(k)$ smallest among all the elements of $\e(k,\x)$, and the order of the computational complexity is typically bounded by $O(M+\tau_0(k)\log{M})$ \cite{Knu73}. Upon obtaining $\a^*_k(\x)$, one only needs to search the minimum of $H(k,\x)$ for $k\in\mathcal{M}$. Therefore, the total searching order is $\sum_{k=1}^M O(M+\tau_0(k)\log{M})$. The worst-case scenario is when $\tau_0(k)\equiv M,\,k=0,1,\dots,M$, which leads to an upper bound of the searching order of $O(M^2\log M)$.\hfill $\square$

\vspace{.20in}
Note that the searching order of $O(M^2\log M)$ is a considerable improvement over $O(2^M)$. This is due to the special form of the loss function and the nature of the parameter space, action space, and multiple decision function space. In a lot of cases, including the specific pairs of loss functions that we will discuss later in Section \ref{sec: concrete}, the searching order may still be lower than $O(M^2\log M)$. 


Observe that in the BMDF described in Theorem \ref{general solution}, we need to obtain the posterior expectation of the form $E_{\btheta|\X=\x}\,[\beta(\btheta\tr\1)\h(\btheta)]$. Recall that $X_m$, given $\btheta$, has the marginal distribution $Q_{m\theta_m}$. We assume for now that the $X_m$s, given $\btheta$, are independent for $m=1,2,\dots,M$. In general, we may also model the $X_m$s to be dependent as will be discussed in Section \ref{sec: data dep}, in which case the computation of the posterior expectation will be discussed in Section \ref{sec: computation}. Denote the density of $Q_{m\theta_m}$ by $q_{m\theta_m}$. Suppose an independent prior distribution of the form described in (\ref{indep prior}) is used. Then, the posterior distribution of $\btheta$ also makes $\theta_m$s independent, with
\begin{equation*}
\pi_m(\theta_m|\x )=\frac{\pi_{m\theta_m}q_{m\theta_m}(x_m)}{\pi_{m0}q_{m0}(x_m)+\pi_{m1}q_{m1}(x_m)} \textrm{  for } m=1,2,\dots,M.
\end{equation*}
Therefore,
\begin{equation*}
E_{\btheta|\x}\,[\beta(\btheta\tr\1)\h(\btheta)]=\sum_{\btheta\in\Theta} \beta(\btheta\tr\1)\h(\btheta) \prod_{m=1}^M \frac{\pi_{m\theta_m}q_{m\theta_m}(x_m)}{\pi_{m0}q_{m0}(x_m)+\pi_{m1}q_{m1}(x_m)}.
\end{equation*}
Notice that, in general, this is a sum of $|\Theta|=2^M$ terms. A particular case is $E(\btheta|\x)=(E(\theta_1|\x),E(\theta_2|\x),\dots,E(\theta_M|\x))$ where, for $m=1,2,\dots,M$,
\begin{equation*}
E(\theta_m|\x)=P(\theta_m=1|\x)=\frac{\pi_{m1}q_{m1}(x_m)}{\pi_{m0}q_{m0}(x_m)+\pi_{m1}q_{m1}(x_m)}.
\end{equation*}

It is worth pointing out that, in general, each component of the posterior expectation $E_{\btheta|\x}\,[\beta(\btheta\tr\1)\h(\btheta)]$ is needed to obtain the BMDF, so that each component of $\delta^*$ may depend on all components of $\X$. This makes the BMDF a {\em{compound}} decision function \cite{SunCai07}. In essence the decision for the $m$th component borrows information from the other components, or as \citeasnoun{Efr10} mentioned, the decision makes use of direct evidence from the $m$th component of the data as well as indirect evidence from other components.

\subsection{BMDF in Composite Hypotheses}\label{subsec: BMDF composite}
Let $\pi$ be a prior density function on the enlarged parameter space $\Theta^\circ$ (see (\ref{extended par vec}) on page \pageref{extended par vec}). Then the Bayes risk of a randomized decision function $\delta\in\mathcal{D}$ is given by $$r_\pi(\delta)=E_{(\btheta,\bgamma,\bxi)} \,E_{\X|(\btheta,\bgamma,\bxi)} \,E_{\A\sim \delta(\X)(\cdot)}\,L(\A,\btheta)=E_\X\, E_{\A\sim \delta(\X)(\cdot)}\, E_{\btheta|\X}\,L(\A,\btheta).$$
Observe that the final form of the Bayes risk is exactly the same as in the simple-vs-simple setting. This implies that the results in Theorem \ref{general solution} apply directly. However, since the parameter space where the prior distribution is defined is now enlarged, the posterior expectation $E_{\btheta|\X}\,L(A,\btheta)$ in the Bayes risk, or $E_{\btheta|\x}\,[\beta(\btheta\tr\1)\h(\btheta)]$ in Theorem \ref{general solution} is now taken with respect to the {\em{marginal}} posterior distribution of $\btheta$, given $\X=\x$.

Assume the $X_m$s, given $\theta$, are independent, for $m=1,2,\dots,M$. Denote the density of $Q_m(\cdot\,;\gamma_m,\xi_m)$, under $H_{m\theta_m}$, by $q_{m\theta_m}(\cdot\,;\gamma_m,\xi_m)$. If an independent prior of the form in (\ref{composite prior}) is assigned, then the marginal posterior distribution of $\btheta$ also makes $\theta_m$s independent, with 
\begin{equation*}
\pi_m(\theta_m|\x) =\frac{\pi_{m\theta_m}\tilde{q}_{m\theta_m}(x_m)}{\pi_{m0}\tilde{q}_{m0}(x_m)+\pi_{m1}\tilde{q}_{m1}(x_m)},
\end{equation*}
where
\begin{equation}\label{tilde q}
\tilde{q}_{m\theta_m}(x_m)=\int_{\Gamma_{m\theta_m}}\int_{\Xi_m}\,p_{m\theta_m}(\gamma_m,\xi_m)q_{m\theta_m}(x_m;\gamma_m,\xi_m){\rm{d}}\xi_m{\rm{d}}\gamma_m.
\end{equation}
Therefore,
\begin{equation*}
E_{\btheta|\x}\,[\beta(\btheta\tr\1)\h(\btheta)]=\sum_{\btheta\in\Theta} \beta(\btheta\tr\1)\h(\btheta) \prod_{m=1}^M \frac{\pi_{m\theta_m}\tilde{q}_{m\theta_m}(x_m)}{\pi_{m0}\tilde{q}_{m0}(x_m)+\pi_{m1}\tilde{q}_{m1}(x_m)}.
\end{equation*}
In particular, $E(\btheta|\x)=(E(\theta_1|\x),E(\theta_2|\x),\dots,E(\theta_M|\x))$ where, for $m=1,2,\dots,M$,
\begin{equation*}
E(\theta_m|\x)=\frac{\pi_{m1}\tilde{q}_{m1}(x_m)}{\pi_{m0}\tilde{q}_{m0}(x_m)+\pi_{m1}\tilde{q}_{m1}(x_m)}.
\end{equation*}
Notice that the integral in $\tilde{q}_{m\theta_m}(x_m)$ may not be in a closed form. Thus Monte Carlo techniques may be needed even in this independent setting. 
Similarly to the simple-vs-simple setting, the BMDF in this composite hypotheses setting is of a compound nature.

\section{Some Concrete Situations}\label{sec: concrete}

\subsection{FP and FN Loss Functions} 
Consider the loss function 
\begin{equation*}
L_{(FP, FN)}(\a,\btheta)=C_0L_{FP}(\a,\btheta)+C_1L_{FN}(\a,\btheta)
             =C_0\frac{\a\tr(\1-\btheta)}{M}+C_1\frac{(\1-\a)\tr\btheta}{M},
\end{equation*}
where $L_{FP}$ and $L_{FN}$ are the false positive proportion and false negative proportion. 
It is clear that the optimal action minimizing $\tilde{L}_{(FP, FN)}$ is
$\a^*_{(FP, FN)}(\x)=((\a^*_{(FP, FN)}(\x))_{m}, m=1,2, \dots, M),$ 
where, for $m=1,2,\dots,M$,
\begin{equation*}
(a^*_{(FP, FN)}(\x))_{m}
       =I\left(\frac{E(\theta_m|\x)}{1-E(\theta_m|\x)}>\frac{C_0}{C_1}\right)
       =I\left(E(\theta_m|\x)>\frac{C_0}{C_0+C_1}\right).
\end{equation*}
The corresponding BMDF $\delta^*$ is such that $\delta^*(\X)=a^*(\X)$. This BMDF is of intuitive form in that the decision on each component is based only on $E(\theta_m|\x)=P(\theta_m=1|\x)$ and the threshold is just $C_0/(C_0+C_1)$. Note that, under $M=1$, this is the Bayes test corresponding to a $C_0/C_1$-loss function \cite{CasBer01}.

\subsection{FDP and FNP Loss Function} 

Consider the loss function 
\begin{equation*}
L_{(FDP,FNP)}(\a,\btheta)=C_0L_{FDP}(\a,\btheta)+C_1L_{FNP}(\a,\btheta)
             =C_0\frac{\a\tr(\1-\btheta)}{(\a\tr\1)\vee 1}+C_1\frac{(\1-\a)\tr\btheta}{((\1-\a)\tr\1)\vee 1},
\end{equation*}
where $L_{FDP}$ and $L_{FNP}$ are the false discovery proportion and the false non-discovery proportion. Note that for this $L\in\mathcal{L}$, $\g_0(\a)=\a$, $\a_0=\1$, $A_1=1$, and $\tau_0(k)=k$.  
Let $(\phi_{(1)}(\x),\phi_{(2)}(\x),\dots,\phi_{(M)}(\x))$ denote the ordered vector associated with $E(\btheta|\x)$. Then,  in Theorem \ref{general solution}, we have
\begin{equation*}
H(k,\x)=\,C_0\dfrac{\sum_{i=1}^k{(1-\phi_{(M-i+1)}(\x))}}{k}+C_1\dfrac{\sum_{i=k+1}^M{\phi_{(M-i+1)}(\x)}}{M-k} .
\end{equation*}
Letting
$k^*(\x)=\argmin_{k\in\mathcal{M}}\,H(k,\x),$
then, by Theorem \ref{general solution}, the optimal action is
\begin{align*}
\a_{(FDP,FNP)}^*(\x)=&(I(r_m(k^*(\x),\x)\leq k^*(\x)),\,m=1,2,\dots,M)\\
        =&(I(\text{rank}(E(\theta_m|\x))\ge M-k^*(\x)+1),m=1,2,\dots,M).
\end{align*}

Notice that for all $k=1,2,\dots,M-1$, $H(k,\x)$ depends only on the ordered vector of $E(\btheta|\x)$, which means that in order to select $k^*(\x)$ we only need to sort $E(\btheta|\x)$ once. Also, the optimal action only requires the rank vector of $E(\btheta|\x)$ after $k^*(\x)$ has been obtained. Therefore, the searching order is reduced to $O(M\log{M})$.  

For this case, the posterior means of the $\theta_m$'s are still the basis of decisions, but in contrast to the previous case, the decision at any particular component depends on {\em{all}} posterior means. One may conclude that the Bayes multiple decision function does not depend on the magnitudes of the $E(\theta_m|\x)$, but rather only on their relative ranks. However, this is not the case since their magnitudes are actually needed to determine $k^*(\x)$.

\subsection{FDP and MDP Loss Functions} 
Consider the loss function 
\begin{equation*}
L_{(FDP, MDP)}(\a,\btheta)=C_0L_{01}(\a,\btheta)+C_1L_{11}(\a,\btheta)
             =C_0\frac{\a\tr(\1-\btheta)}{(\a\tr\1)\vee 1}+C_1\frac{(\1-\a)\tr\btheta}{(\btheta\tr\1)\vee 1},
\end{equation*}
where $L_{01}$ and $L_{11}$ are the false discovery proportion and the missed discovery proportion. Same as for the pair of FDP and FNP, $\g_0(\a)=\a$, $\a_0=\1$, $A_1=1$, and $\tau_0(k)=k$. 
For $k=1,2,\dots,M$, let
\begin{equation*}
\tilde{e}(k,\x)=\dfrac{C_0}{k}E(\btheta|\x)+C_1E\left( \frac{\btheta}{(\btheta\tr\1)\vee \1}|\x\right),
\end{equation*}
and denote by $(\tilde{e}_{(1)}(k,\x),\tilde{e}_{(2)}(k,\x),\dots,\tilde{e}_{(M)}(k,\x))$ the ordered vector of $\tilde{e}(k,\x)$.
Then
\begin{equation*}
H(k,\x)=\begin{cases}
         C_1\1\tr E\left(\frac{\btheta}{(\btheta\tr\1)\vee \1}|\x\right) &\text{if $k=0$}\\
         C_1\1\tr E\left(\frac{\btheta}{(\btheta\tr\1)\vee \1}|\x\right)+C_0-\sum_{i=1}^k \tilde{e}_{(M-i+1)}(k,\x) &\text{if $k\neq 0$}
        \end{cases},
\end{equation*}
and
\begin{equation*}
k^*(\x)=\begin{cases}
        0 &\text{if $C_0>\max_{k\in\mathcal{M}\setminus  \{ {0} \}}\sum_{i=1}^k \tilde{e}_{(M-i+1)}(k,\x)$} \\
        \argmax_{k\in\mathcal{M}\setminus  \{ {0} \}}\sum_{i=1}^k \tilde{e}_{(M-i+1)}(k,\x) &\text{otherwise}
        \end{cases}.
\end{equation*}
By Theorem \ref{general solution}, the optimal Bayes action is
\begin{equation*}
\a^*(\x)=(I(\text{rank}(\tilde{e}_m(k^*(\x),\x))\ge M-k^*(\x)+1),\, m=1,2,\dots,M).
\end{equation*}
Observe that $a^*(\x)$ and $k^*(\x)$ depend on the values and ranks of $\tilde{e}_m(k,\x)$ for $k\in\mathcal{M}$. The searching order in this case is $O(M^2\log{M})$. It can be shown that when the posterior probability distribution of $\btheta$ specifies independent components, the searching order is reduced to $O(M\log{M})$.

\section{Dependent Data Structure}\label{sec: data dep}

In section \ref{sec: BMDF}, formulas of the posterior expectations are given under the assumption that the $X_m$'s, $m=1,2,\dots,M$, are independent. However in various situations this assumption is far from being realistic. Therefore, we consider a dependent data structure, namely in the general setting introduced in section \ref{general problem}, the distribution of $\X$ given $(\btheta,\gamma_\btheta)$, denoted by $Q_\btheta(\cdot;\gamma_\btheta)$ with $\bvartheta(Q_\btheta(\cdot;\gamma_\btheta))=\btheta$.


In the simple hypotheses setting, the goal is to specify $Q_\btheta$, the joint probability distribution of $\X$ given $\btheta$, such that the marginal probability distribution of $X_m$ is $Q_m=Q_{m\theta_m}$. Let $\mathcal{M}_0(\btheta)=\{m\in\{1,2,\dots,M\}\colon\theta_m=0\}$ and $\mathcal{M}_1(\btheta)=\{m\in\{1,2,\dots,M\}\colon\theta_m=1\}$. Assume that $X_{\mathcal{M}_0(\btheta)}\equiv \{X_m\colon m\in\mathcal{M}_0\}$ is a collection of independent random vectors, and $X_{\mathcal{M}_1(\btheta)}\equiv \{X_m\colon m\in\mathcal{M}_1\}$ is a collection of possibly dependent vectors, and the collections $X_{\mathcal{M}_0(\btheta)}$ and  $X_{\mathcal{M}_1(\btheta)}$ are independent of each other. Moreover, borrowing survival analysis ideas, we assume that the dependence structure of the collection $X_{\mathcal{M}_1(\btheta)}$ is that induced by a frailty-based model. Specifically, assume that
$$Q_\btheta\left(\left.\bigcap_{m\in\mathcal{M}_1(\btheta)} [X_m\in B_m]\,\right| Z=z \right)=\prod_{m\in\mathcal{M}_1(\btheta)}[\breve{Q}_m(B_m)]^z,$$
for all $B_m\in\mathcal{B}_m$, $m\in\mathcal{M}_1(\btheta)$, where $\breve{Q}_m$'s are some distributions on $\mathcal{X}_m$'s, and the frailty $Z\in\mathcal{Z}$ is assumed to have a distribution $G$. Therefore, the joint distribution of $X_{\mathcal{M}_1(\btheta)}$ is
\begin{equation}\label{dep joint}
Q_\btheta\left(\bigcap_{m\in\mathcal{M}_1(\btheta)} [X_m\in B_m]\right)=\int_\mathcal{Z}\prod_{m\in\mathcal{M}_1(\btheta)}[\breve{Q}_m(B_m)]^zG({\rm{d}}z),
\end{equation}
for all $B_m\in\mathcal{B}_m$, $m\in\mathcal{M}_1(\btheta)$. Recall that in the simple-vs-simple multiple hypotheses testing setting, under $H_{m1}$, $X_m\sim Q_{m1}$ marginally. So the distributions $\breve{Q}_m$, $m\in\mathcal{M}_1(\btheta)$, should be such that these conditions are satisfied. Let $\mathcal{L}_G$ be the Laplace transform of the distribution function $G$, that is, $\mathcal{L}_G(u)=\int_\mathcal{Z} e^{-uz}G({\rm{d}}z),\, \forall u\in\mathbb{R}.$ Denote $M_1\equiv|{\mathcal{M}_1(\btheta)}|$. The following result gives the joint distribution of the dependent collection $X_{\mathcal{M}_1(\btheta)}$ in terms of the collection of the marginal distributions $\{Q_{m1}:m\in\mathcal{M}_1(\btheta)\}$ under the assumed frailty-based model.

\begin{prop}\label{prop:copula}
The frailty-based model described in (\ref{dep joint}) is an $M_1$-dimensional Archimedean copula $C_G$ such that $$Q_\btheta\left(\bigcap_{m\in\mathcal{M}_1(\btheta)} [X_m\in B_m]\right)=C_G(Q_{m1}(B_m),m\in\mathcal{M}_1(\btheta))$$ for all $B_m\in\mathcal{B}_m$, where $C_G\colon [0,1]^{M_1}\to[0,1]$ is defined via
$$C_G(u_1,u_2,\dots,u_{M_1})=\mathcal{L}_G\left(\sum_{m=1}^{M_1}\mathcal{L}_G^{-1}(u_m)\right).$$
\end{prop}

{\noindent\bf Proof: } According to the model, marginally for $m\in\mathcal{M}_1(\btheta)$ and all $B_m\in\mathcal{B}_m$,
\begin{align*}
Q_{m1}(B_m)&=\int_\mathcal{Z} \breve{Q}_m(B_m)^zG({\rm{d}}z)\\
&=\int_\mathcal{Z}\exp{\big\{-z\big\{-\log{\breve{Q}_m(B_m)}\big\}\big\}} G({\rm{d}}z)=\mathcal{L}_G(-\log\breve{Q}_m(B_m)).
\end{align*}
Thus $\breve{Q}_m(B_m)=\exp(-\mathcal{L}_G^{-1}(Q_{m1}(B_m))$. So,

\begin{align*}
Q_\btheta\left(\bigcap_{m\in\mathcal{M}_1(\btheta)} [X_m\in B_m]\right)&=\int_\mathcal{Z}\prod_{m\in\mathcal{M}_1(\btheta)}(\exp(-\mathcal{L}_G^{-1}(Q_{m1}(B_m))))^zG({\rm{d}}z)\\
    &=\int_\mathcal{Z}\exp\left(-\left[\sum_{m\in\mathcal{M}_1(\btheta)}\mathcal{L}_G^{-1}(Q_{m1}(B_m))\right]\cdot z\right)G({\rm{d}}z)\\
    &=\mathcal{L}_G\left(\sum_{m\in\mathcal{M}_1(\btheta)} \mathcal{L}_G^{-1}(Q_{m1}(B_m))\right).
\end{align*}

To show $C_G$ is an Archimedean Copula, it is sufficient to show that the function $\mathcal{L}_G^{-1}$ is a \textit{strict generator} of a copula, which is that it is a continuous strictly decreasing convex function from $[0,1]$ to $[0,\infty]$ with $\mathcal{L}_G^{-1}(1)=0$ and $\mathcal{L}_G^{-1}(0)=\infty$ \cite{Nel99}. But these are straight-forward to verify using properties of the Laplace transform. \hfill $\square$\\

\vspace{-.25in}
Thus, a frailty-induced dependent full data model is given by
\begin{equation}\label{dep model}
Q_\btheta\left(\bigcap_{m=1}^M [X_m\in B_m]\right)=\left(\prod_{m\in\mathcal{M}_0(\btheta)}Q_{m0}(B_m)\right) C_G\Big[Q_{m1}(B_m),m\in\mathcal{M}_1(\btheta)\Big],
\end{equation}
for all $B_m\in\mathcal{B}_m$, $m=1,2,\dots,M$. 
Notice that the distribution function $G$ may have nuisance parameters, say, $G(\cdot)=G(\cdot,\bupsilon)$, where $\bupsilon\in\Upsilon$. In this case, to calculate the posterior expectations, a prior on $\Upsilon$ will also be needed.

In the composite hypotheses testing setting, we are to specify $Q_\btheta(\cdot;\bgamma,\bxi)$, the joint probability distribution of $\X$ given $(\btheta,\bgamma,\bxi)$, such that the marginal probability distribution of $X_m$ is $Q_m=Q_{m\theta_m}(\cdot;\gamma_{m},\xi_m)$. The result in Proposition \ref{prop:copula} is easily extended to get
$$Q_\btheta\left(\bigcap_{m\in\mathcal{M}_1(\btheta)} [X_m\in B_m];\bgamma,\bxi\right)=C_G\Big[Q_{m1}(B_m;\gamma_{m},\xi_m),m\in\mathcal{M}_1(\btheta)\Big],$$ and the full data model is therefore given by 
$$Q_\btheta\left(\bigcap_{m=1}^M [X_m\in B_m];\bgamma,\bxi\right)=\left(\prod_{m\in\mathcal{M}_0(\btheta)}Q_{m0}(B_m;\gamma_{m},\xi_m)\right) C_G\Big[Q_{m1}(B_m;\gamma_{m},\xi_m),m\in\mathcal{M}_1(\btheta)\Big].$$

\section{Sequential Monte Carlo}\label{sec: computation}
The applicability of the algorithm of finding the Bayes optimal action in Theorem \ref{general solution} is contingent on an efficient way of calculating the posterior expectations 
$$E(H(\btheta)|\X=\x)=\sum_{\btheta\in\Theta}H(\btheta)\pi(\btheta|\x),$$
where $H(\btheta)$ takes the form $\beta(\btheta\tr\1)h(\btheta)$ and $\pi(\btheta|\x)$ is the posterior probability mass function of $\btheta$, given data $\x$. For example, the specific forms of the $H$ function desired in case of FDP and MDP loss functions are 
$H(\btheta)=\btheta \textrm{ and } H(\btheta)=\frac{\btheta}{(\btheta\tr\1)\vee 1}.$  
As pointed out in Sections \ref{sec: BMDF} and \ref{sec: data dep} Monte Carlo integration is needed for approximating the posterior expectations. However, in the regular Importance Sampling (IS) algorithms \cite{Rip87}, as $M$, the dimension of $\btheta$, increases, the computational complexity of calculating the weights also increases. So it is important to consider a sequential application of the importance sampling methods \cite{Gor93}. Notice that $m=1,2,\dots,M$ does not necessarily represent time or positions in an ordered sequence, and that the proposed dependent data structure is not necessarily a state space model as what is usually the case in sequential importance sampling applications, but the technique provides a visual solution through $m=1,2,\dots,M$ that deals with the dimensionality and monitors the efficiency of the sampling procedure. 

\subsection{Simple Hypotheses}

Let $\pi$ be a prior probability mass function of $\btheta$. Under the dependent data model described in (\ref{dep model}), the desired posterior expectation is given by
$$I(\x)=E(H(\btheta)|\x)=\frac{\sum_{\btheta\in\Theta}H(\btheta)\pi(\btheta)Q_\btheta({\rm{d}}\x)}{\sum_{\btheta\in\Theta}\pi(\btheta)Q_\btheta({\rm{d}}\x)},$$
where $Q_\btheta({\rm{d}}\x)=\left(\prod_{m\in\mathcal{M}_0(\btheta)}q_{m0}(x_m)\right)C_G(q_{m1}(x_m), m\in\mathcal{M}_1(\btheta))$, and $q_{m0}$ and $q_{m1}$ are the density functions of $Q_{m0}$ and $Q_{m1}$, respectively. Consider an independent data-adaptive trial probability mass function $g$ of $\btheta$ given by $g(\btheta)=g(\btheta|\x)=\prod_{m=1}^M g_m(\theta_m|x_m),$ where, for $m=1,2,\dots,M$, $g_m(\theta_m|x_m)$ is the marginal posterior of $\theta_m$, given $X_m=x_m$. This is a Bernoulli distribution of form
\begin{equation}\label{simp trial}
g_m(\theta_m|x_m)\propto\pi_m(\theta_m)q_{m0}(x_m)^{1-\theta_m}q_{m1}(x_m)^{\theta_m},
\end{equation}
where $\pi_m(\theta_m)$ is the marginal prior probability of $\theta_m$.
Denote by $\x_{1:m}=(x_1,x_2,\dots,x_m)$ and $\btheta_{1:m}=(\theta_1,\theta_2,\dots,\theta_m)$ for $m=1,2,\dots,M$. The joint prior distribution of $\btheta$ can be written as $\pi(\btheta)=\prod_{m=1}^M \pi_m(\theta_m|\btheta_{1:m-1})$, where $\pi_1(\theta_1|\btheta_{1:0})$ is the marginal distribution of $\theta_1$, and for $m=2,3,\dots,M$, $ \pi_m(\theta_m|\btheta_{1:m-1})$ is the probability of $\theta_m$ given $\btheta_{1:m-1}$ under the joint prior $\pi$. Notice that if $\pi$ is independent, then $ \pi_m(\theta_m|\btheta_{1:m-1})=\pi_m(\theta_m).$ However, dependent prior structures can also be constructed by the frailty-induced models similarly to the dependent data structure, in which case $ \pi_m(\theta_m|\btheta_{1:m-1})$ is not necessarily reduced to $\pi_m(\theta_m).$ 

Let $q_{\btheta_{1:m}}(\x_{1:m})$ be the marginal density function of $\X_{1:m}$, given $\theta_{1:m}$, under the dependent data structure proposed in (\ref{dep model}). Since the dependent structure is induced by a frailty model, we have, for $m=2,3,\dots,M$, $q_{\btheta_{1:m}}(\x_{1:m})=q_{\btheta_{1:m}}(x_m|\x_{1:m-1})\cdot q_{\btheta_{1:m-1}}(\x_{1:m-1})$, where 
\begin{equation*}
q_{\btheta_{1:m}}(x_m|\x_{1:m-1})=\frac{q_{\btheta_{1:m}}(\x_{1:m})}{q_{\btheta_{1:m-1}}(\x_{1:m-1})}=
\begin{cases}
q_{m0}(x_m) &\text{if $\theta_m=0$}\\
\dfrac{C_G(q_{t1}(x_t), t\in\mathcal{M}_1(\btheta_{1:m}))}{C_G(q_{t1}(x_t), t\in\mathcal{M}_1(\btheta_{1:m-1}))}&\text{if $\theta_m=1$}
\end{cases}.
\end{equation*}
Thus the full density of the data has the sequential form  
$q_\btheta(\x)=\prod_{m=1}^M q_{\btheta_{1:m}}(x_m|\x_{1:m-1})$
with $q_{\theta_1}(x_1|\x_{1:0})=1$. So there is a recursive formula for calculating the importance weight function, given by 
$w_m(\btheta_{1:m}|\x_{1:m})=w_{m-1}(\btheta_{1:m-1}|\x_{1:m-1})u_m(\btheta_{1:m}|\x_{1:m})$, 
where the increment 
$u_m(\btheta_{1:m}|\x_{1:m})=\dfrac{q_{\btheta_{1:m}}(x_m|\x_{1:m-1})\pi(\theta_m|\btheta_{1:m-1})}{g_m(\theta_m|x_m)}$
satisfies
\begin{equation}\label{simp increment}
u_m(\btheta_{1:m}|\x_{1:m})\propto 
\begin{cases}
\dfrac{\pi(\theta_m=0|\btheta_{1:m-1})}{\pi_m(\theta_m=0)} &\text{if $\theta_m=0$}\\
\dfrac{C_G(q_{t1}(x_t), t\in\mathcal{M}_1(\btheta_{1:m}))}{C_G(q_{t1}(x_t), t\in\mathcal{M}_1(\btheta_{1:m-1}))\cdot q_{m1}(x_m)} \cdot  \dfrac{\pi(\theta_m=1|\btheta_{1:m-1})}{\pi_m(\theta_m=1)} &\text{if $\theta_m=1$}
\end{cases}.
\end{equation}
Note that $w_M(\btheta_{1:M}|\x_{1:M})=w(\theta|\x)$, the importance weight function. In summary, we can now sample particles and calculate the importance weights in a sequential manner.

The fundamental difficulty of SIS is the degeneracy of the weights. For large values of $M$, the weights $w^{(r)}(\btheta|\x)$, $r=1,2,\dots,R$, are all close to $0$ except for one that is close to $1$, which would result in a poor estimate eventually. One of the solutions to this difficulty is through resampling, or using the so-called the bootstrap filter \cite{Liu01}, in which, after resampling, all the importance weights are set to $1/R$ so that we make sure all particles make important contributions to the MC estimate. But resampling at each $m=1,2,\dots,M$ may be computationally expensive, so we would only resample whenever the empirical effective sample size is too low. Via this procedure, we can make sure that the weights will not diverge.

\begin{alg}{\bf{(SMC in Simple Hypotheses) }}
\begin{enumerate}
\item Fix a large integer $R$ and a threshold $\rho\in (0,1]$.
\item Iterate for $m=1,2,\dots,M$.
\begin{enumerate}
\item For all $r=1,2,\dots,R$, generate $\theta_m^{(r)}$ independently from $g_m(\cdot|x_m)$ in (\ref{simp trial}), and set $\btheta_{1:m}^{(r)}=(\btheta_{1:m-1}^{(r)},\theta_m^{(r)})$.
\item Compute the increments $u_m(\btheta_{1:m}^{(r)}|\x_{1:m})$ given in (\ref{simp increment}), and the importance weights $w_m^{(r)}\equiv w_m(\btheta_{1:m}^{(r)}|\x_{1:m})=w_{m-1}(\btheta_{1:m-1}^{(r)}|\x_{1:m-1})u_m(\btheta_{1:m}^{(r)}|\x_{1:m})$.
\item Compute $ESS=R/ \sum_{r=1}^R(w_m^{(r)})^2$.
\item If $ESS<\rho R$, normalize the weights, and resample, with replacement, $R$ particles from $\{\btheta_{1:m}^{(r)}\colon r=1,2,\dots,R\}$ according to the normalized weights, and set all the weights to $1/R$.
\end{enumerate}
\item $I(\x)=E(H(\btheta)|\x)$ is approximated by $\hat{I}(\x)=\dfrac{\sum_{r=1}^R H(\btheta_{1:M}^{(r)})w_M^{(r)}}{\sum_{r=1}^R w_M^{(r)}}$.
\end{enumerate}
\end{alg}

\subsection{Composite Hypotheses}
Let $\pi$ be the prior probability function on the enlarged parameter space $\Theta^\circ$ described in (\ref{composite prior}), and denote the independent Bernoulli marginal prior probability on $\Theta$ by $\pi_0$. Then, under the dependent data model described in (\ref{dep model}), the desired posterior expectation is given by
$$I(\x)=E(H(\btheta)|\x)=\frac{\sum_{\btheta\in\Theta}H(\btheta)\pi_0(\btheta)\tilde{Q}_\btheta({\rm{d}}\x)}{\sum_{\btheta\in\Theta}\pi_0(\btheta)\tilde{Q}_\btheta({\rm{d}}\x)},$$
where $\tilde{Q}_\btheta({\rm{d}}\x)=\left(\prod_{m\in\mathcal{M}_0(\btheta)}\tilde{q}_{m0}(x_m)\right) C_G\Big[\tilde{q}_{m1}(x_m), m\in\mathcal{M}_1(\btheta)\Big]$, and $\tilde{q}_{m0}$ and $\tilde{q}_{m1}$ are as defined in (\ref{tilde q}). 
Consider a data-adaptive trial density $\tilde{g}$ on the enlarged parameter space, given by 
$\tilde{g}(\btheta,\bgamma,\bxi)=\tilde{g}(\btheta,\bgamma,\bxi|\x)=\prod_{m=1}^M\tilde{g}_m(\theta_m,\gamma_m|x_m),$ where, for $m=1,2,\dots,M$,
$$\tilde{g}_m(\theta_m,\gamma_m,\xi_m|x_m)\propto (\pi_{m0}\breve{q}_{m0}(x_m)p_{m0}(\gamma_m,\xi_m))^{1-\theta_m}(\pi_{m1}\breve{q}_{m1}(x_m)p_{m1}(\gamma_m,\xi_m))^{\theta_m},$$
and for $j=0,1$, $\breve{q}_{mj}(x_m)=q_{mj}(x_m;\hat{\gamma}_{mj}(x_m),\hat{\xi}_{mj}(x_m))$, where $\hat{\gamma}_{mj}(x_m)$ and $\hat{\xi}_{mj}(x_m)$ are some convenient estimates, for example, maximum likelihood or method-of-moments estimates, of $\gamma_m$ and $\xi_m$, under the marginal models $X_m\sim Q_{mj}$. 
Notice that if $\breve{q}$ is replaced by $\tilde{q}$, $\tilde{g}_m$ would be equal to the marginal posterior of $(\theta_m,\gamma_m,\xi_m)$ given $x_m$. Since $\breve{q}_{mj}(x_m)$ is an approximation of $\tilde{q}_{mj}(x_m)$, $\tilde{g}_m$ also provides some guidance to the posterior distribution by making use of the data.
Write the prior in (\ref{composite prior}) as $\pi(\btheta,\bgamma,\bxi)=\prod_{m=1}^M \pi_m(\theta_m,\gamma_m,\xi_m)$, where for $m=1,2,\dots,M$, $\pi_m(\theta_m,\gamma_m,\xi_m)=[\pi_{m0}p_{m0}(\gamma_m,\xi_m)]^{1-\theta_m}[\pi_{m1}p_{m1}(\gamma_m,\xi_m)]^{\theta_m}$. Similarly to the simple-vs-simple case, the full density of data has a sequential form given by 
$q_\btheta(\x; \bgamma, \bxi)=\prod_{m=1}^M q_{\btheta_{1:m}}(x_m|\x_{1:m-1}; \bgamma_{1:m},\bxi_{1:m}),$ where 
\begin{equation*}
q_{\btheta_{1:m}}(x_m|\x_{1:m-1}; \bgamma_{1:m},\bxi_{1:m})=
\begin{cases}
q_{m0}(x_m; \gamma_m,\xi_m) &\text{if $\theta_m=0$}\\
\dfrac{C_G(q_{t1}(x_t; \gamma_t,\xi_t), t\in\mathcal{M}_1(\btheta_{1:m}))}{C_G(q_{t1}(x_t; \gamma_t,\xi_t), t\in\mathcal{M}_1(\btheta_{1:m-1}))}&\text{if $\theta_m=1$}
\end{cases}.
\end{equation*}
So the recursive formula for calculating the importance weight function is 
$$w_m(\btheta_{1:m}|\x_{1:m}; \bgamma_{1:m},\bxi_{1:m})=w_{m-1}(\btheta_{1:m-1}|\x_{1:m-1}; \bgamma_{1:m-1},\bxi_{1:m-1})u_m(\btheta_{1:m}|\x_{1:m}; \bgamma_{1:m},\bxi_{1:m}),$$
where the increment satisfies
\begin{equation}\label{comp increment}
u_m(\btheta_{1:m}|\x_{1:m}; \bgamma_{1:m-1},\bxi_{1:m-1})\propto 
\begin{cases}
\dfrac{q_{m0}(x_m;\gamma_m,\xi_m)}{\breve{q}_{m0}(x_m)} &\text{if $\theta_m=0$}\\
\dfrac{C_G(q_{t1}(x_t; \gamma_t,\xi_t), t\in\mathcal{M}_1(\btheta_{1:m}))}{C_G(q_{t1}(x_t; \gamma_t,\xi_t), t\in\mathcal{M}_1(\btheta_{1:m-1}))\cdot \breve{q}_{m1}(x_m)}&\text{if $\theta_m=1$}
\end{cases}.
\end{equation}
Therefore, the SIS algorithm is very similar to that for the simple-vs-simple case with the increment replaced by (\ref{comp increment}).

\section{Illustration and Simulations}\label{sec: illustration}



\subsection{Composite Alternatives with Independent Gaussian Observations}
Assume $X_i|\theta_i$ are independent $N(\mu_{\theta_i},1)$, where $\mu_{0i}=0$ and $\mu_{1i}\stackrel{iid}{\sim} N(0,\sigma^2)$ for fixed $\sigma$. Consider tests $H_{0i}:\mu_i=0$ vs $H_{1i}:\mu_i\neq 0$ with independent prior $(\theta_i,\mu_i)\sim (1-\pi)I(\theta_i=0)+\pi\phi(\mu_i;0,\sigma^2)I(\theta_i=1)$, where $\pi$ is the fixed prior probability of the alternative hypotheses, and $\phi(\cdot;\mu,\sigma^2)$ is the density function of a normal distribution with mean $\mu$ and variance $\sigma^2$. 
Then the posterior means are
$$E(\theta_i|X_i)=\frac{\pi \phi(X_i;0,1+\sigma^2)}{(1-\pi) \phi(X_i;0,1)+\pi \phi(X_i;0,1+\sigma^2)}.$$
We performed 1000 simulations with $M=12$, true $\sigma=4$ and true proportion of alternatives $\pi=0.5$ for all three procedures. For both correct prior parameters with $\pi^*=0.5$, $\sigma^*=4$ and misspecified prior parameters with $\pi^*=0.7$, $\sigma^*=10$, the empirical FDP, FNP and MDP were calculated for $C_0/C_1=1,2$. The results in Table \ref{Simu Com} compares the BMDF associated with the three pairs of  loss functions to the Benjamini and Hochberg (BH) procedure \cite{BenHoc95}  with a false discovery rate (FDR) threshold 0.05. To find the BMDF associated with the (FDP, MDP) pair of loss functions, the posterior expectation $E(\btheta/(\btheta\tr\1\vee 1)|\x)$ was calculated {\em{exactly}} by using a recursive formula and {\em{approximately}} by a Monte Carlo approximation. Note that the use of BH procedure in the simulation study is to enable comparison with the most commonly-used, albeit frequentist, multiple testing procedure. However, it worth mentioning that it is not totally fair to compare the BMDF to the BH procedure because they are designed under different criteria. 
In Table \ref{Simu Com}, the empirical risks of the three BMDF are comparable to those of the BH procedure when the cost ratio $C_0/C_1=2$, and with a smaller cost ratio ($=1$), the empirical FDP become larger and the empirical FNP and MDP become smaller since with a smaller cost the BMDF sacrifices larger Type I error to achieve an optimal combined risks. With misspecified prior parameters, the combined empirical risks stay almost the same as those with the correctly specified prior parameters, but the Type I error FDP becomes smaller and the Type II error FNP and MDP become larger. This is because we specified a higher prior probability of the alternatives than the truth, which results in more discoveries. Finally, the exactly and approximately-calculated posterior expectations result in similar BMDF associated with the (FDP, MDP) pair of loss functions in terms of the empirical risks. Since the exact calculation is too computationally expensive for large $M$, in the following illustration we will utilize the Monte Carlo approximated posterior expectations. 

 \begin{table}[h]
\begin{center}
\begin{tabular}{|c|c|c||c|c|c|c|c||c|} \hline
$\pi^*$ & $\sigma^*$ & $\lambda$ & Loss      &   (FP, FN)    &  (FDP, FNP)   & \multicolumn{2}{c||}{(FDP,MDP)} & BH \\ 
& & & Functions & & & \multicolumn{1}{c}{Exact} & Approx. & \\\hline
\multirow{6}{*}{0.5} & \multirow{6}{*}{4} & & $\hat{FDP}$ &  0.115 &  0.052  &  0.184 & 0.112  & 0.025\\ 
& & 1 & $\hat{FNP}$  &  0.260  & 0.278   &  0.211 & 0.253 &  0.309 \\ 
& & & $\hat{MDP}$  &  0.323  &  0.422  & 0.230  &  0.311 &  0.450 \\ \cline{3-9} 
& &  & $\hat{FDP}$ &  0.051 &  0.010 &  0.035 & 0.026 & 0.025\\ 
& & 2 & $\hat{FNP}$ &  0.290 &   0.324 &  0.297 & 0.303 & 0.309 \\ 
& & & $\hat{MDP}$ &  0.388 &  0.484  & 0.412  &  0.428 &  0.450\\ \hline\hline
\multirow{6}{*}{0.7} & \multirow{6}{*}{10} & & $\hat{FDP}$ &  0.098 & 0.049   & 0.165  & 0.100 &0.027\\ 
& &1 & $\hat{FNP}$ &  0.270 & 0.287   &  0.225 & 0.263  & 0.304\\ 
& & & $\hat{MDP}$  & 0.332  &   0.424 &  0.248 &  0.324 &0.448 \\\cline{3-9} 
& &  & $\hat{FDP}$ & 0.063  &  0.017  &  0.042 & 0.035 & 0.027\\ 
& & 2 & $\hat{FNP}$ &  0.285 &  0.317  &  0.290 & 0.297 & 0.304\\ 
& & & $\hat{MDP}$ & 0.390  &  0.483  &  0.410 & 0.428 & 0.448\\ \hline
\end{tabular}
\caption[Average of empirical FDP, FNP and MDP in 1000 simulations in a composite multiple testing problem with independent Gaussian observations.]{Average of empirical FDP, FNP and MDP in 1000 simulations in a composite multiple testing problem with independent Gaussian observations. True prior parameters are $\pi=0.5$ and $\sigma=4$.}
\label{Simu Com}
\end{center}
\end{table}

\subsection{Simple-vs-Simple with Dependent Exponential Observations}
Consider a situation where, for all $m$, $X_{m1}, X_{m2}, \ldots, X_{mn}$ are IID with 
$\mathcal{Q}_{m0} = \{EXP(\lambda_m): \lambda_m=\lambda_0\}$, and
$\mathcal{Q}_{m1} =  \{EXP(\lambda_m): \lambda_m=\lambda_1\}$. 
A Gamma$(\kappa,\kappa)$ frailty induces dependency among $\{\X_m:m\in\mathcal{M}_1\}$.
The independent Bernoulli prior with $P(\theta_m=1)=\pi$ is used. Data was generated under the true model parameters:
$M = 500, n=30, \pi = .30, \lambda_0 = 1, \lambda_1= .5, \kappa = 2 $. 
To find the BMDF, a frailty-model with exponential marginals and a Gamma($\kappa,\kappa$) frailty is used.
In the Sequential Monte Carlo, R = 1000 particles are used with $\pi^* = .20, \kappa^* = 3.$
Results in one replicate for different pairs of loss functions is shown in Table \ref{Illu Exp}. The cost ratio used is
$C_0/C_1 = 1.$ Notice that the performances of all three BMDF are satisfactory even under misspecifications of the prior probabilities of the alternative hypotheses and the hyperparameter in the distribution of the frailty. 

\begin{table}[h]
\center
\begin{tabular}{||c||c|c||c|c||c|c||} \hline
Loss Functions Pair & \multicolumn{2}{c||}{(FP, FN)} & \multicolumn{2}{c||}{(FDP,FNP)} &  \multicolumn{2}{c||}{(FDP,MDP)} \\ \hline
Actions & Nulls & Alts & Nulls & Alts & Nulls & Alts \\ \hline
0: Accepts & 334 & 4 & 331 & 4 & 338 & 6 \\
1: Rejects & 6 & 156 & 9 & 156 & 2 & 154 \\ \hline
\end{tabular}
\caption{Results in one simulation in an exponential multiple testing problem.}
\label{Illu Exp}
\end{table}

\subsection{Two Group Composite Hypothesis with Independent Gaussian Observations}\label{subsec: two group simul}
\begin{center}
\begin{figure}[h!]
\includegraphics[height=3.1in,width=6.7in,angle=0]{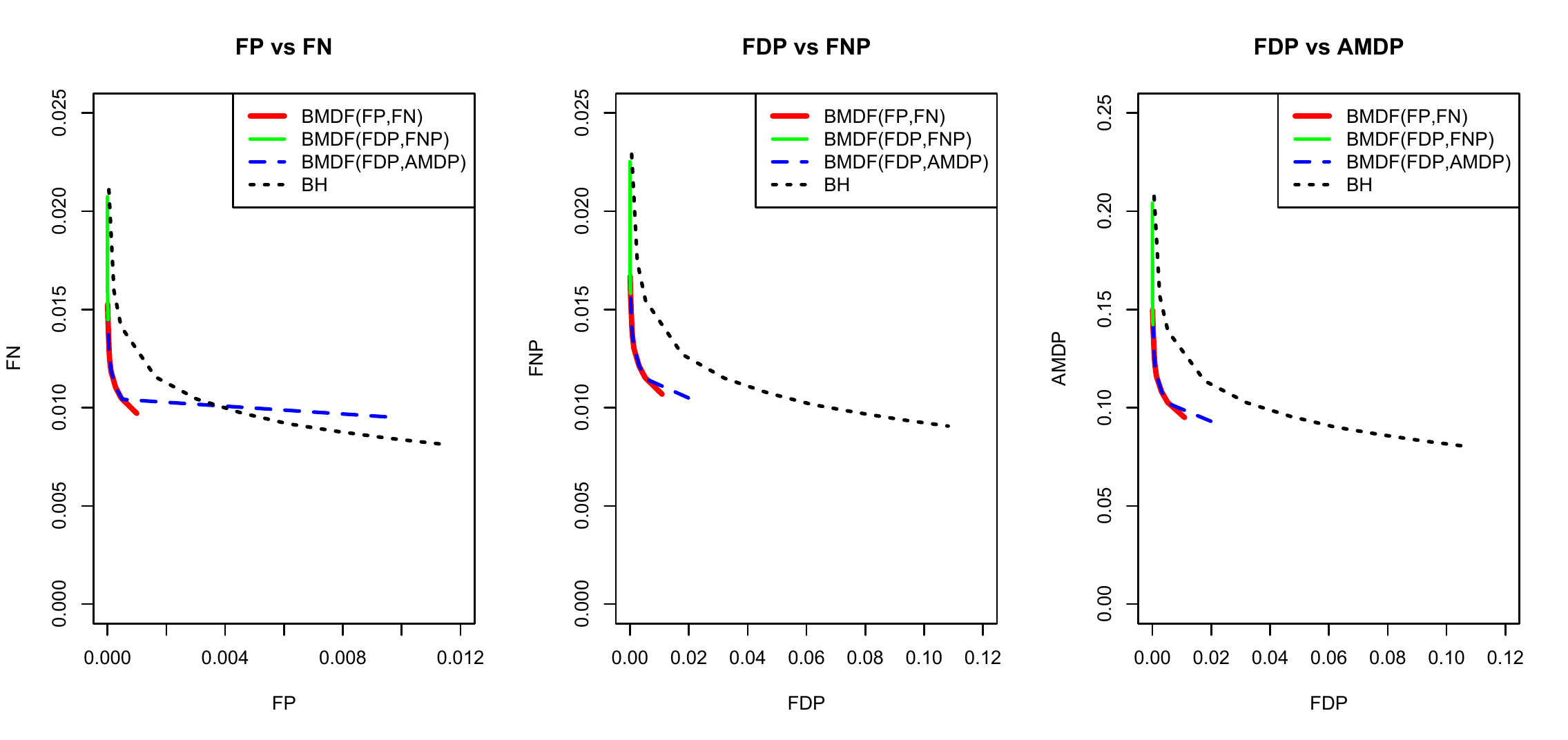}
\caption[Graphs of empirical risk pairs for three BMDF with different loss functions when cost ratio varies and the BH procedure when the FDR threshold varies in a two-group composite hypothesis problem with independent Gaussian observations and correct prior parameters are specified.]{Graphs of ($\hat{FP}$, $\hat{FN}$), ($\hat{FDP}$, $\hat{FNP}$), and ($\hat{FDP}$, $\hat{AMDP}$) for three BMDF with different loss functions when cost ratio varies and the BH procedure when the FDR threshold varies in a two-group composite hypothesis problem with independent Gaussian observations. Correct prior parameters are specified. }
\label{Illu Twogroup Correctprior}
\end{figure}
\end{center}

\begin{center}
\begin{figure}[h!]
\includegraphics[height=3.1in,width=6.7in,angle=0]{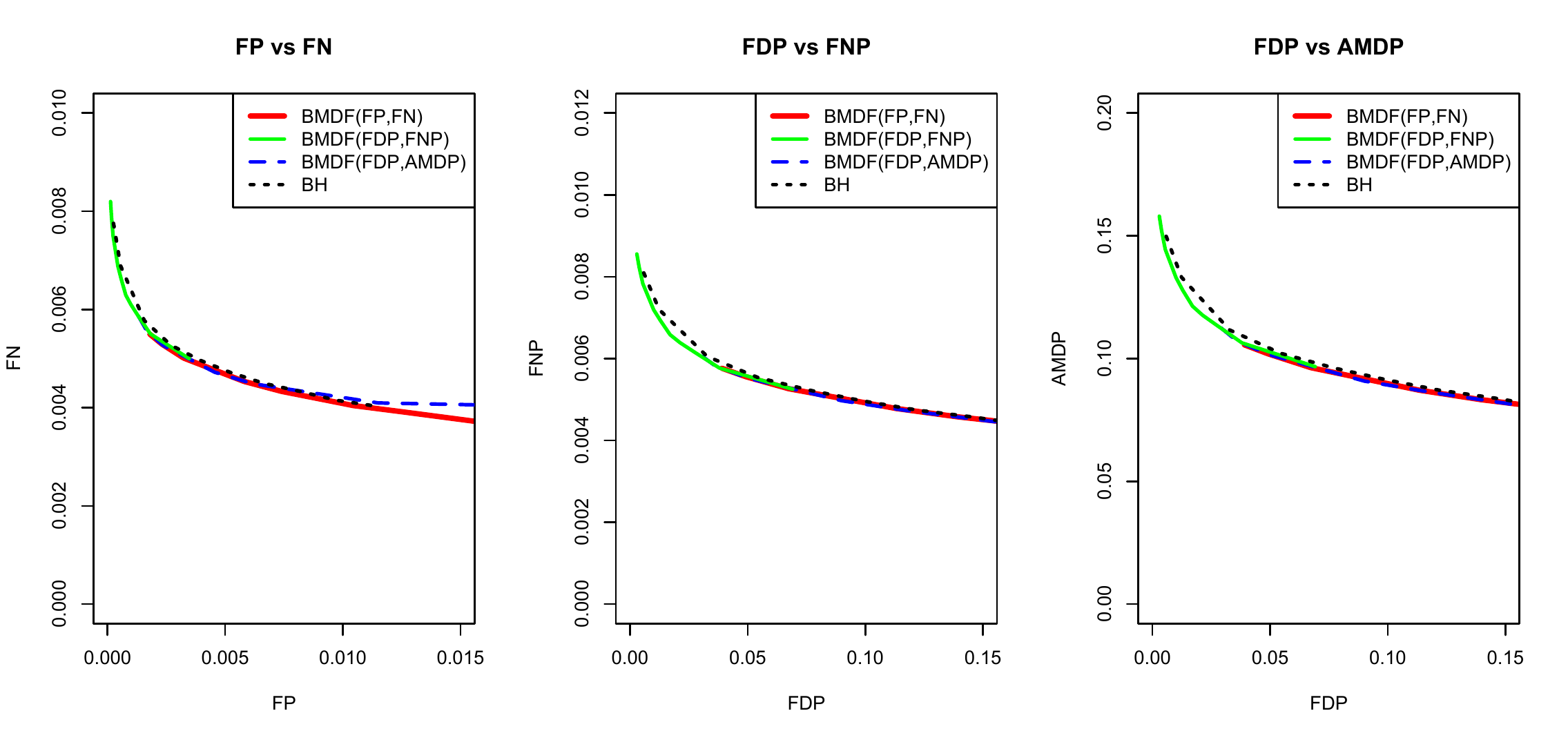}
\caption[Graphs of empirical risk pairs for three BMDF with different loss functions when cost ratio varies and the BH procedure when the FDR threshold varies in a two-group composite hypothesis problem with independent Gaussian observations and the prior parameters are partially misspecified and partially empirically estimated..]{Graphs of ($\hat{FP}$, $\hat{FN}$), ($\hat{FDP}$, $\hat{FNP}$), and ($\hat{FDP}$, $\hat{AMDP}$) for three BMDF with different loss functions when cost ratio varies and the BH procedure when the FDR threshold varies in a two-group composite hypothesis problem with independent Gaussian observations. The prior parameters are partially misspecified and partially empirically estimated.}
\label{Illu Twogroup EB}
\end{figure}
\end{center}

Assume that
$\mathbf{X}_m = (X_{m1}, X_{m2}, \ldots, X_{mn_1}, Y_{m1}, Y_{m2}, \ldots, Y_{mn_2})$
are independent with $X_{mi}\stackrel{iid}{\sim}N(\mu_{m1},\sigma_m^2)$ and $Y_{mi}\stackrel{iid}{\sim}N(\mu_{m2},\sigma_m^2)$. Consider tests $H_{0m}: \mu_{m1} = \mu_{m2}$ vs $H_{1m}:\mu_{m1} \neq \mu_{m2}$ with independent Bernoulli prior on $\btheta$ with probability $\pi$ for the alternatives,  and conjugate prior for nuisance parameters given by, 
$\mbox{for } m\in\mathcal{M}_1$, $\mu_{m1},\mu_{m2}\stackrel{iid}{\sim}N(\nu_m, k_0\sigma_m^2)$, $\sigma_m^{-2}\sim \mbox{Gamma}(\alpha, \beta)$, 
and $\mbox{for } m\in\mathcal{M}_0$, $\mu_{m1}=\mu_{m2}\stackrel{iid}{\sim} N(\nu_m, k_0\sigma_m^2)$, $\sigma_m^{-2}\sim \mbox{Gamma}(\alpha, \beta)$. True parameters used to generate data are $M=500, \pi=0.1, k_0=200, \alpha=4, \beta=4, \nu=20. $ We performed 1000 simulations each with correct prior parameters and misspecified prior parameters: $ \pi^*=0.05, k_0^*=100, \alpha^*=20$ as well as {\em{empirically}} estimated other parameters via $\nu={\overline{\x}_{m}}, \beta=S^2(\x_m)(\alpha-1)/(k_0+1))$.  In order to stabilize the results when $\pi$ is small, we used an adjusted version of  MDP (AMDP) given by $L(\a,\btheta)=\dfrac{(\1-\a)\tr\btheta}{(\btheta\tr\1)+ 1}$. We implemented the three BMDF associated with different loss functions with 10 different cost ratios $C_0/C_1$ and BH procedure with 10 different FDR thresholds. Figure \ref{Illu Twogroup Correctprior} and Figure \ref{Illu Twogroup EB} show graphs of three empirical risk pairs for all four procedures with correctly specified prior parameters and partially misspecified and partially empirically-estimated prior parameters, respectively. In Figure  \ref{Illu Twogroup Correctprior} with the correct prior specification, the empirical risk curves for all three BMDF are well below that for the BH procedure indicating better performance, except for the BMDF associated with (FDP, AMDP) pair of loss functions where the average FP loss is surprisingly high. This is because with small prior probability $\pi$ for the alternatives, the AMDP could be large for some simulated data, and  when the cost ratio $C_0/C_1$ is very small, the BMDF would rather sacrifice a large number of false positives to achieve the optimal combined risk of FDP and AMDP. In Figure \ref{Illu Twogroup EB}, when the prior parameters are partially misspecified and partially empirically estimated, the results of all four procedures almost coincide. However, more study about the empirical Bayes ideas is needed in future research. 

\subsection{A Real Microarray Data Analysis}

\begin{center}
\begin{figure}[h!]
\includegraphics[height=5.3in,width=6.7in,angle=0]{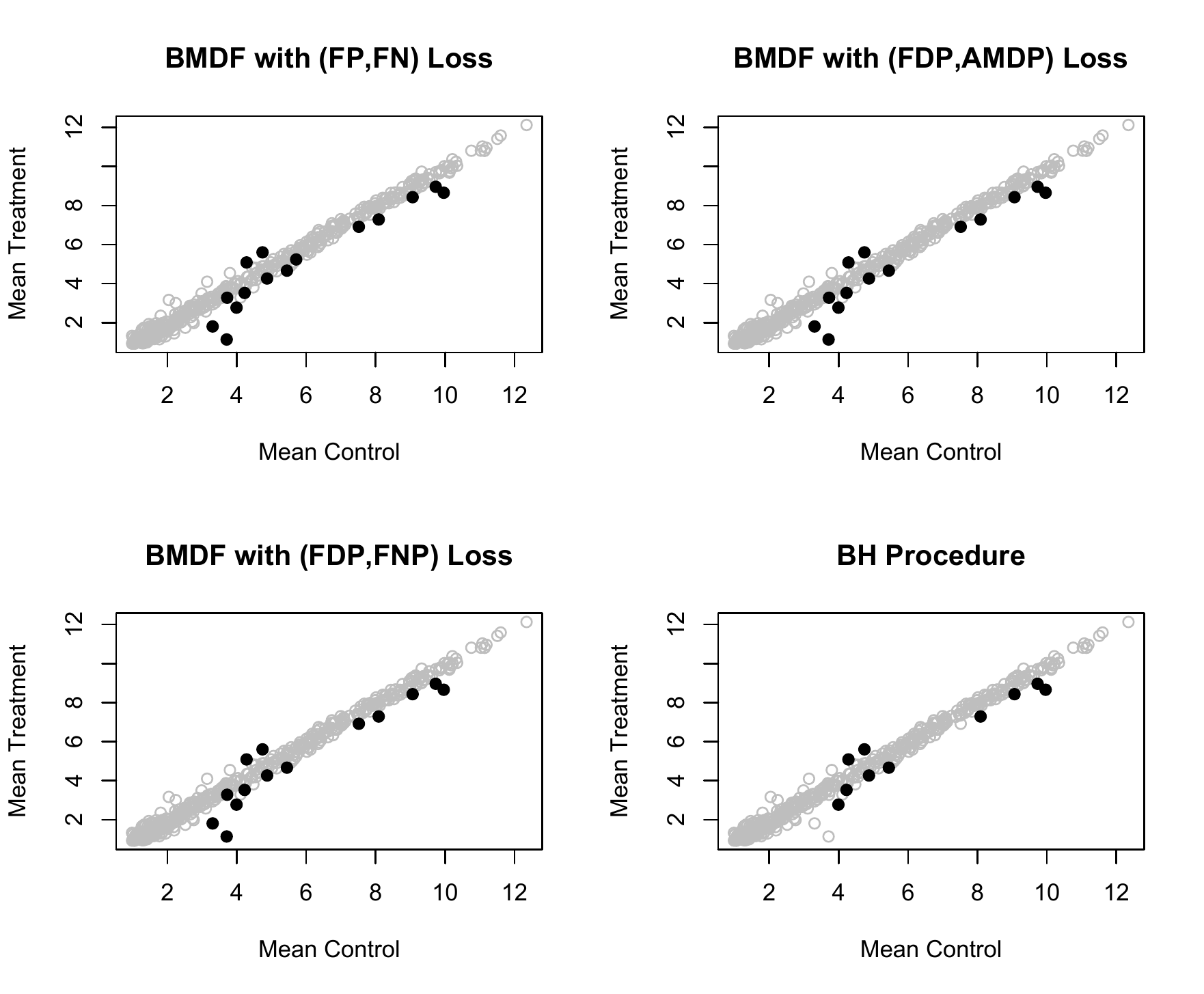}
\caption[Mean expression level of the control group versus the treatment group for 500 genes in a microarray data set. ]{Mean expression level of the control group versus the treatment group for 500 genes in a microarray data set. Gray circles are non discoveries and black circles are discoveries. Top left: BMDF with (FP, FN) loss function pair; top right: BMDF with (FDP, AMDP) loss function pair; bottom left: BMDF with (FDP,FNP) loss function pair; bottom right: BH procedure.}
\label{Illu CellPlot}
\end{figure}
\end{center}

In a colon cancer tumor metastasis study conducted in the laboratory of Dr. Marge Pe\~na at the University of South Carolina, expression levels for 41268 genes from mice tissues were obtained through a Agilent Technology microarray. For each gene, five replicates were obtained for a control group and five replicates for a metastatic group.
For our illustration, we randomly selected 500 genes out of the 41268 on which to apply our BMDF. We assume the independent two-group Gaussian model, and use the partially empirically-estimated prior parameters described in Section \ref{subsec: two group simul}. Cost ratios used for (FP, FN), (FDP, FNP) and (FDP, AMDP) loss function pairs are 
$C_0/C_1 = 3, 0.2, 2$ respectively, and for BH procedure the FDR threshold is 0.05. These cost ratios and thresholds were chosen according to the simulation results in Section \ref{subsec: two group simul} in order for the four procedures to have similar empirical FDP. Out of the 500 genes, the BMDF associated with (FP, FN), (FDP, FNP) and (FDP, AMDP) loss function pairs found 15, 14 and 14 genes differentially expressed across groups, respectively, and the BH procedure found 10. Figure \ref{Illu CellPlot} shows the mean expression level of the control group versus the treatment group of the discovered genes. Notice that all 10 genes discovered by the BH procedure were also discovered by the three BMDFs.

\section{Concluding Remarks}\label{sec: conclude}
BMDF defines a class of multiple decision procedures that achieve optimality in the Bayesian framework associated with a general class of loss functions. The results in Theorem \ref{general solution} describe the form of the BMDF and provide an efficient algorithm of finding it in multiple testing settings. Notice that the pairs of loss function are not limited to the ones described in Section \ref{sec: concrete}. For example, the adjusted MDP given by $L(\a,\btheta)=\dfrac{(\1-\a)\tr\btheta}{(\btheta\tr\1)+ 1}$ may help stabilize the Bayes optimal actions when the prior probabilities of the alternative hypotheses $\pi_m=Pr(\theta_m=1)$ are small. Also notice that  the choice of the loss function pairs and the cost ratio should be pre-determined. 

The frailty-based dependent data structure is a class of flexible models if the distribution of the frailty is modeled in a hierarchical manner with hyperparameters. Similarly, the prior distribution could also be dependent with frailty-based structures. Furthermore, the SMC could be easily implemented in the computations of the posterior expectations.  However, not all dependent structures is frailty-based. Therefore, in real data analysis, model validation is needed to see the validity of the dependent structure.

One possible extension of the research in this paper is in two-class prediction problems where the form of BMDF could be extended with the loss functions replaced by prediction errors. Other future studies include the extension to model selection problems, and also in the empirical Bayes approach for determining prior hyperparameter values. In particular, it is of interest to study whether the empirical Bayes procedures are equivalent in some sense to the non-Bayesian BH multiple decision functions and the procedure in \citeasnoun{PenHabWu10}.

\section{Acknowledgements}
The authors thank Dr.\,Marge Pe\~na for providing the dataset and discussing microarrays with us. The authors acknowledge support from National Science Foundation (NSF) Grant
DMS 0805809, National Institutes of Health (NIH) Grant RR17698, and Environmental
Protection Agency (EPA) Grant RD-83241902-0 to the University of Arizona with subaward
number Y481344 to the University of South Carolina. This paper is based on a portion of the first author's PhD dissertation. 


\bibliographystyle{ECA_jasa}
\bibliography{MHTrefs}

\end{document}